\documentclass[12pt]{article}

\usepackage{amsmath,amsfonts,latexsym}
\makeatletter
\newdimen\mainfontsize \mainfontsize=1\@ptsize pt
\topskip=\mainfontsize
   \@maxdepth=\maxdepth

\makeatother

\begin{document}

\newtheorem{theorem}{Theorem}
\newtheorem{cor}[theorem]{Corollary}
\newtheorem{lemma}[theorem]{Lemma}
\newtheorem{prop}[theorem]{Proposition}
\newcommand{\bfs}{\mbox{${\bf s}$}}
\newcommand{\cP}{\mbox{$\mathcal P$}}
\newcommand{\cC}{\mbox{$\mathcal C$}}
\newcommand{\T}{\mbox{$\mathcal T$}}
\newcommand{\bft}{\mbox{$\bf t$}}
\def\enpf{
   {  \parfillskip=0pt\hfil {\hbox{$\Box$}} \par\bigskip  }
   }

\def\R{{\mathbb R}}
\def\C{{\mathbb C}}
\def\Z{{\mathbb Z}}
\def\P{{\mathbb P}}
\def\E{{\mathbb E}}
\def\A{{\cal A}}
\def\F{{\cal F}}
\def\ef{ \tilde{\cal B}^{\theta,b}}
\def\i-c{\otimes}
\def\s-c{\odot}

\begin{center}
\begin{Large}
\begin{bf}
\hbox{Random matrices, non-colliding processes and queues}
\end{bf}
\end{Large}

\bigskip

{\sc Neil O'Connell}

Laboratoire de Probabilit\'es, Paris 6

175 rue du Chevaleret, 75013 Paris

{\tt neil@ccr.jussieu.fr}

(To appear in {\em S\'eminaire de Probabilit\'es XXXVI})
\end{center}

\begin{quote}
This is survey of some recent results connecting random matrices,
non-colliding processes and queues.
\end{quote}

\section{Introduction}

It was recently discovered by Baik, Deift and Johansson~\cite{bdj}
that the asymptotic distribution of the length of the longest increasing
subsequence
in a permutation chosen uniformly at random from $S_n$, properly
centred and normalised, is the same as the asymptotic distribution
of the largest eigenvalue of an $n\times n$ GUE random matrix, properly centred
and normalised, as $n\to\infty$.  This distribution
had earlier been identified by Tracy and Widom~\cite{tw} in the random
matrix context, and it is now known as the Tracy-Widom law.

The length of the longest increasing subsequence in a random permutation
can be thought of as a `last-passage time' for a certain directed percolation
problem on the plane;  this directed percolation problem is closely
related to the one-dimensional totally asymmetric exclusion process
(with low density of particles travelling at high speed) or equivalently,
an infinite series of $M/M/1$ queues in tandem (with low density of
customers and high service rates).  On the other hand, the eigenvalues
of a GUE random matrix of dimension $n$ have the same law as the positions,
after a unit length of time, of $n$ independent standard Brownian motions
started from the origin and conditioned (in the sense of Doob) never to
collide.  These remarks are for the purpose of convincing the reader
that there might be some interesting connections between random matrices,
non-colliding processes and queues.  Indeed there are, and that is
the topic of this paper.

Let us concentrate on the following, more exact, connection between directed
percolation and random matrices which was more recently observed by
Baryshnikov~\cite{bar} and Gravner, Tracy and Widom~\cite{gtw}.
Let  $B=(B_1,\ldots,B_n)$ be a standard $n$-dimensional Brownian motion,
and write $B(s,t)=B(t)-B(s)$ for $s<t$.
\begin{theorem}
\label{bgtw}
{\rm [Baryshnikov; Gravner-Tracy-Widom]}
The random variable
\begin{equation}\label{M_n}
M_n = \sup_{0\le s_1\le\cdots\le s_{n-1}\le 1}
\left\{ B_n(0,s_1)+\cdots +B_1(s_{n-1},1)\right\} ,
\end{equation}
has the same law as the largest eigenvalue of an $n$-dimensional GUE random
matrix.
\end{theorem}
The proofs of this fact given in~\cite{bar} and~\cite{gtw} are based on the
Robinson-Shensted-Knuth (RSK) correspondence,
and do not make use of the facts that
\begin{itemize}
\item[{\rm (a)}] $M_n$ has a queueing interpretation~\cite{gw}, and
\item[{\rm (b)}] the largest eigenvalue of a GUE random matrix has an
interpretation in terms
of non-colliding Brownian motions~\cite{dyson,grabiner}.
\end{itemize}
In~\cite{oy2}, a proof is given of a
more general result which is based entirely on these interpretations.
This more general result identifies a path-transformation $\Gamma_n(B)$
of $B$, which has the same law as that of $n$ independent Brownian motions
conditioned (in the sense of Doob) never to collide.  This process, which
we denote by $\hat B$, is the eigenvalue process associated with Hermitian
Brownian
motion and $\hat B(1)$ has the same distribution as the eigenvalues of an
$n$-dimensional GUE random matrix~\cite{dyson,grabiner}.
The transformed process $\Gamma_n(B)$ has the property that
its largest component at time $1$ is given by $M_n$, and so
Theorem~\ref{bgtw} follows.

At the heart of the proof, which will be outlined in this paper, is a
celebrated theorem of classical queueing theory, which states that,
in equilibrium, the output of a stable M/M/1 queue is Poisson.
This result is usually attributed to Burke, who gave the first proof in 1956.
In what follows, we shall also refer to it as an `output theorem'.
It follows from the reversibility of the M/M/1 queue.
There is an easy extension of this theorem,
which can be proved by similar reversibility arguments; when phrased
in the setting of `max-plus algebra', this extension immediately yields
the two-dimensional result.  The result in higher dimensions is then obtained
by considering a series of queues in tandem and applying an induction
argument.  The statement of Theorem~\ref{bgtw} seems considerably less
mysterious
(to me at least!)
in the setting of queues and non-colliding processes.

In the case $n=2$, the fact that $\Gamma_n(B)$ and $\hat B$ have the
same law is equivalent to Pitman's representation for the three-dimensional
Bessel process.  In the case $n=3$, it yields a representation for planar
Brownian motion conditioned to stay forever in a wedge of angle $\pi/3$;
a partial converse of this result was discovered earlier by Biane~\cite{biane}.

As we remarked above, the process $\hat B$ has the same law as the
eigenvalue process
associated
with Hermitian Brownian motion.  Bougerol and Jeulin~\cite{bj} have
(independently)
obtained a similar representation for this process, by completely different
methods,  which is also consistent with Theorem~\ref{bgtw}.
Their results are presented in the more general setting of
Brownian motion on symmetric spaces.
The relationship between these two
representations will be discussed elsewhere.

There are certain symmetries in the max-plus algebra which play a crucial role,
and these symmetries also exist in other algebras, including the conventional
algebra on the reals.  As a consequence, there are analogues of these output
theorems in many different settings, and it seems that there are some general
phenomena at work.  In all cases, symmetry and reversibility play a key role.
We will present some of these examples and their implications.

We will also describe briefly how output theorems and the study of tandem
systems
can be used to obtain first order asymptotic results for various directed
percolation and directed polymer models.

The outline of the paper is as follows.  In the next section we present
some background material on random matrices and non-colliding Brownian motions.
In section~\ref{state} we give a precise statement of the representation
theorem obtained in~\cite{oy2}, and briefly discuss the special cases
$n=2$ and $n=3$.  In section~\ref{burke} we describe the output
theorem for the M/M/1 queue, how it is proved using reversibility, and how
it can be extended using similar reversibility arguments.
In section~\ref{ncp} we show how the extended version of Burke's
theorem can be used to obtain a representation for non-colliding Poisson
processes, and how the representation for non-colliding Brownian motions
follows.
In section~\ref{output}
we discuss output theorems generally and give some examples.
In section~\ref{perc} we describe how output theorems can be
used to  obtain first order asymptotic results for various directed
percolation and directed polymer models.

\section{Random matrices and non-colliding Brownian motions}

Recall that a $n\times n$ GUE random matrix $A\in\C^{n\times n}$ is constructed
as follows:  it is Hermitian, that is, $A=A^*(=(\bar A)^t)$;
the entries $\{ A_{ij},\ i\le j\}$ are independent;
on the diagonal $A_{ii}$ are standard real normal random variables;
below the diagonal, $\{ A_{ij},\ i< j\}$ are standard complex
normal random variables, that is, the real and imaginary parts of
$A_{ij}$ are independent centered real normal random variables,
each with variance $1/2$;  above the diagonal we set $A_{ji}=\bar A_{ij}$.
Here, $\bar z=x-iy$ denotes the complex conjugate of $z=x+iy$.
The joint density (with respect to Lebesgue measure on $\R^n$)
of the eigenvalues of $A$ is given by Weyl's formula
(see, for example,~\cite{mehta}):
\begin{equation}\label{weyl}
f_{GUE}(\lambda_1,\ldots,\lambda_n)=Z^{-1}\prod_{i\ne j}(\lambda_i-\lambda_j)
\prod_i\exp(-\lambda_i^2/2),
\end{equation}
where $Z=(2\pi)^{n/2}\prod_{j=1}^{n-1}j!$.

This distribution has an interpretation in terms of non-colliding Brownian
motions.
The (Vandermonde) function $h$ defined by
\begin{equation}\label{h}
h(x)=\prod_{i<j}(x_j-x_i)
\end{equation}
is harmonic on $\R^n$, and moreover, is a strictly positive harmonic
function for Brownian motion killed when it exits the Weyl chamber
\begin{equation}
W=\{ x\in\R^n:\ x_1<x_2<\cdots <x_n\} .
\end{equation}
For $x\in \R^n$, let $\P_x$ denote the law of $B$ started at $x$;
for $x\in W$, let $\hat{\P}_x$ denote the law of the $h$-transform of $B$
started at $x$, where
$h$ is given by~(\ref{h}).  The laws $\hat{\P}_x$ and $\P_x$ are related
as follows.  If $T$ denotes the first exit time of $B$ from $W$,
and $\F_t$ the natural filtration of $B$, then for $A\in\F_t$,
\begin{equation}
\hat{\P}_x (A) = \E_x \left( \frac{h(B_t)}{h(x)} 1_{A\cap\{T>t\}} \right) .
\end{equation}
It is well-known, and easy to check using the Karlin-MacGregor formula
(see, for example,~\cite{ko}) that
\begin{equation}\label{j}
\lim_{W\ni x\to 0} \hat{\P}_x (X_t\in dy) = C_t h(y)^2 \P_0(X_t\in dy),
\end{equation}
where
\begin{equation}\label{ct}
C_t=\left[ t^{n(n-1)/2} \prod_{j=1}^{n-1} j! \right]^{-1}
\end{equation}
is a normalisation constant.  Note that the RHS of (\ref{j}) is equal
to $f_{GUE}(y)$ when $t=1$.

There is a more sophisticated connection at the process level.
Note that we can define the law $\hat{\P}_{0+}$,
for $A\in \T_{t}=\sigma(B_u,\ u\ge t)$, $t>0$, by
\begin{equation}
\hat{\P}_{0+} (A) = \P_0 \left[ C_t h(B_t)^2 \hat{\P}_{B_t} (\theta_t A)
\right] ,
\end{equation}
where $\theta$ is the shift operator (so that $\theta_t A \in \T_0$).
The fact that this is well-defined follows from~(\ref{j}).
Now, {\em Hermitian Brownian motion} is constructed in the same way as a GUE
random matrix, but with Brownian motions instead of normal random variables.
It is a fact that the law of process of eigenvalues of Hermitian Brownian
motion
is given by $\hat{\P}_{0+}$ (see, for example,~\cite{dyson,grabiner,mehta}).

For related work on non-colliding diffusions and random matrices,
see~\cite{biane,cl,hobson,ko}, and references therein.

\section{A generalisation of Theorem~\ref{bgtw}}\label{state}

Let $D_0(\R_+)$ denote the space of cadlag paths $f:\R_+\to\R$
with $f(0)=0$.  For $f,g\in D_0(\R_+)$,
define $f\i-c g\in D_0(\R_+)$ and $f\s-c g\in D_0(\R_+)$ by
\begin{equation}\label{ic}
(f\i-c g)(t)=\inf_{0\le s\le t}[f(s)+g(t)-g(s)] ,
\end{equation}
and
\begin{equation}\label{sc}
(f\s-c g)(t)=\sup_{0\le s\le t}[f(s)+g(t)-g(s)] .
\end{equation}
Unless otherwise deleniated by parentheses, the default order of operations is
from left to right; for example, when we write $f\i-c g\i-c h$, we mean
$(f\i-c g)\i-c h$.
Define a sequence of mappings $\Gamma_k:D_0(\R_+)^k\to D_0(\R_+)^k$
by
\begin{equation}\label{g2}
\Gamma_2(f,g)=(f\i-c g, g\s-c f),
\end{equation}
and, for $k>2$,
\begin{align}\label{gn}
\Gamma_k(f_1, &  \ldots,f_k)=(f_1\i-c f_2\i-c\cdots\i-c f_k,\nonumber\\
& \Gamma_{k-1}(f_2\s-c f_1,f_3\s-c (f_1\i-c f_2),\ldots,f_k\s-c (f_1\i-c
\cdots\i-c f_{k-1}))) .
\end{align}
The operations $\i-c$ (and $\s-c$) arise naturally in a queueing context;
they can also be regarded as operator products
in the min-plus (resp. max-plus) algebra.  We refer the
reader to~\cite{bb} for more about max-plus algebra and its applications.

Let $\hat B$ be a realisation of $\hat \P_{0+}$, as defined in the previous
section.
Then the largest eigenvalue of a $n\times n$ GUE random matrix has the same law
as $\hat B_n(1)$, and Theorem~\ref{bgtw} states that $M_n$ and
$\hat B_n(1)$ have the same law.

Now observe that, if $\Gamma_n(B)_1$ denotes the first component in
the $n$-dimensional process $\Gamma_n(B)$, then
\begin{eqnarray}
\Gamma_n(B)_1(t) &=& (B_1\i-c\cdots\i-c B_n)(t) \\
&=& \inf_{0\le s_1\le\cdots\le s_{n-1}\le t}
\left\{ B^{(1)}_{(0,s_1)}+\cdots +B^{(n)}_{(s_{n-1},t)}\right\} ;
\end{eqnarray}
thus, by symmetry, Theorem~\ref{bgtw} states that $\Gamma_n(B)_1(1)$
has the same law as the smallest eigenvalue of an $n$-dimensional
GUE random matrix or, equivalently, the random variable $\hat B_1(1)$.
We remark that, although it is not immediately obvious, it can also be
shown that
\begin{equation}
\Gamma_n(B)_n = B_n\s-c\cdots\s-c B_1 ,
\end{equation}
from which it also follows that $M_n=\Gamma_n(B)_n(1)$ has the same law
as $\hat B_n(1)$.  In fact~\cite{oy2}:
\begin{theorem}\label{main}
The processes $\Gamma_n(B)$ and $\hat B$ have the same law.
\end{theorem}
The proof of this identity in law, given in~\cite{oy2}, is based
on an extension of the output theorem for
a stable $M/M/1$ queue in equilibrium, and certain related symmetries in the
max-plus algebra.

In the next section we discuss the output theorem for M/M/1 queues.
In Section~\ref{ncp} we show how one can deduce the analogue of
Theorem~\ref{main} for non-colliding Poisson processes.  This is
first done in the case of unequal rates, where the non-collision
conditioning is non-singular.  The case of equal rates follows,
and it is then purely a technical exercise to apply Donsker's theorem
and recover the Brownian version.  We remark that the Poisson case
with equal rates is interesting in its own right, as it is
closely related to the Charlier orthogonal polynomial ensemble.

As we remarked in the introduction, the case $n=2$ is equivalent
to Pitman's representation for the three-dimensional Bessel process,
which states that, if $X$ is a standard one-dimensional Brownian motion
and $M(t)=\max_{0<s<t}X(t)$, then $2M-X$ is a three-dimensional Bessel process.
The three-dimensional Bessel process, being the $h$-transform of Brownian
motion on $\R_+$ with $h(x)=x$, can also be interpreted as Brownian
motion conditioned to stay positive.
To see the connection, first note that in this case, $R=(\hat B_2-\hat
B_1)/\sqrt{2}$
is a three-dimensional Bessel process and $(\hat B_2+\hat B_1)/\sqrt{2}$
is a standard Brownian motion, independent of $R$.  On the other hand,
$$\Gamma_2(f,g)_2+\Gamma_2(f,g)_1=g+f,$$
and
$$\Gamma_2(f,g)_2-\Gamma_2(f,g)_1=g\s-c f - f\i-c g = 2m-x,$$
where $x=g-f$ and $m(t)=\max_{0<s<t}x(s)$.

We now consider the statement of Theorem~\ref{main} in the case $n=3$.
In this case Brownian motion in the Weyl chamber
$$W_3=\{x\in\R^3:\ x_1<x_2<x_3\}$$
can be mapped onto
Brownian motion in a cone of angle $\pi/3$ on the plane.  Let
$$C=\{ z\in\C:\ |\arg z| < \pi/6\},$$
and define $\phi:W_3\to C$ by $$\phi(x_1,x_2,x_3)=\left(
\frac{x_3-x_1}{\sqrt{2}} ,
\frac{2x_2-x_1-x_3}{2\sqrt{3}} \right) .$$
If $B$ is a standard Brownian motion in $\R^3$, then $\phi(B)$ is a standard
Brownian motion in $\C$, and moroever, $\phi(W_3)=C$ (see, for
example,~\cite{biane,ou}).
It therefore follows from Theorem~\ref{main} that $\phi(\Gamma_3(B))$
is a Brownian motion conditioned (in the sense of Doob) to stay in $C$ forever.
As remarked by Biane~\cite{biane}, there is only one way to condition on
this event.

In fact, Biane~\cite{biane} proved that, if $X$ is a standard Brownian
motion in $\C$
conditioned to stay in $C$ forever, and $J_t=\inf_{u>t}X(u)$, then $\Im(3J-X)$
is a standard one-dimensional Brownian motion.  Thus, if we set
$X=\phi(\Gamma_3(B))$,
then $\Im(3J-X)$ is a standard one-dimensional Brownian motion.  The expression
$\Im(3J-X)$ can be simplified considerably, but I do not see how to write it as
a linear function of $B$; thus, it is not clear at this point how to recover
Biane's result directly from Theorem~\ref{main}.

\section{An extension of Burke's theorem}\label{burke}

The stationary M/M/1 queue can be constructed as follows.
Let $A$ and $S$ be independent Poisson processes on $\R$
with respective intensities $0<\lambda <\mu$.  For intervals $I$,
open, half-open or closed, we will denote by $A(I)$ the measure of $I$
with respect to $dA$; for $I=(0,t]$ we will simply write $A(t)$, with
the convention that $A(0)=0$.  Similarly for $S$ and any other point process
we introduce.
For $t\in\R$, set
\begin{equation}\label{past}
Q(t) = \sup_{s\le t} [A(s,t]-S(s,t]]^+,
\end{equation}
and for $s<t$,
\begin{equation}\label{dmass}
D(s,t]=A(s,t]+Q(s)-Q(t).
\end{equation}
In the language of queueing theory, $A$ is the {\em arrivals} process,
$S$ is the {\em service} process, $Q$ is the {\em queue-length} process,
and $D$ is the {\em departure} process.  With this construction it is also
natural (and indeed very important for what follows) to define the {\em
unused service} process by
\begin{equation}
U=S-D .
\end{equation}

We will use the following notation for reversed processes.  For a point process
$X$, the reversed process $\bar X$ is defined by $\bar{X}(s,t)=X(-t,-s).$
The reversed queue-length process $\bar Q$ is defined to be the
right-continuous modification of $\{ Q(-t),\ t\in\R\}$.

Burke's theorem states that $D$ is a homogeneous Poisson process
with intensity $\lambda$.  On a historical note, this fact was anticipated
by O'Brien~\cite{obrien} and Morse~\cite{morse}, and proved in 1956
by Burke~\cite{burke}.   Independently, it was also proved by
Cohen~\cite{cohen}.
In 1957, Reich~\cite{reich} gave the following
very elegant proof which uses reversibility.  The process
$Q$ is reversible (in fact, all stationary birth and death processes are
reversible).
It follows that the joint law of $A$ and $D$ is the same as the
joint law of $\bar D$ and $\bar A$. In particular, $\bar D$, and
hence $D$, is a Poisson process with intensity $\lambda$.

Burke also proved that, for each $t$, $\{ D(s,t],\ s\le t\}$
is independent of $Q(t)$.  This property is now called {\em
quasi-reversibility}.  Note that it also follows from Reich's
reversibility argument.  Discussions on Burke's theorem and related
material
can be found in the books of Br\'emaud~\cite{bremaud,bremaud2},
Kelly~\cite{kelly} and Robert~\cite{robert}.

Set
\begin{equation}\label{tmass}
T=A+U.
\end{equation}
\begin{theorem}\label{basic}
The point processes $D$ and $T$ are independent Poisson processes with
respective intensities $\lambda$ and $\mu$.
\end{theorem}
{\bf Proof.}
First note that, given $Q$, $U$ is a homogeneous Poisson process
with intensity $\mu$ on the set $I=\{s\in\R:\ Q(s)=0\}$, and if we let
$V$ be another Poisson process with intensity $\mu$ on the complement
of $I$, which is conditionally independent of $U$ given $Q$, then
(unconditionally) $N=U+V$ is a homogeneous Poisson process
with intensity $\mu$ on $\R$ which is independent of $Q$.
Now, $(A,S)$ can be written as a simple function of $(Q,N)$,
$(A,S)=\varphi(Q,N)$ say.  By construction, we have $(\bar D,\bar T)
=\varphi(\bar Q,\bar N)$.  Now we use the reversibility of $Q$ and
$N$ to deduce that $(\bar D,\bar T)$, and hence $(D,T)$, has the
same law as $(A,S)$, as required.
\hfill $\Box$

Note that
\begin{equation}\label{future}
Q(t)=\sup_{u>t}[D(t,u)-T(t,u)].
\end{equation}
We also have, on $\{Q(0)=0\}$,
\begin{equation}
\{ (D(t),T(t)),\ t\ge 0\} = \{ \Gamma_2(A,S)(t),\ t\ge 0\} .
\end{equation}

Theorem~\ref{basic} has the following multi-dimensional extension,
which relates to a sequence of M/M/1 queues in tandem.
Let $A,S_1,\ldots,S_n$ be independent Poisson processes
with respective intensities $\lambda, \mu_1, \ldots, \mu_n$, and assume
that $\lambda < \min_{i\le n} \mu_i$.
Set $D_0=A$ and, for $k\ge 1$, $t\in\R$, set
\begin{equation}\label{pastk}
Q_k(t) = \sup_{s\le t} [D_{k-1}(s,t]-S_k(s,t]]^+,
\end{equation}
and for $s<t$,
\begin{equation}\label{dmassk}
D_k(s,t]=D_{k-1}(s,t]+Q_k(s)-Q_k(t),
\end{equation}
\begin{equation}
T_k(s,t]=S_k(s,t]-Q_k(s)+Q_k(t).
\end{equation}
\begin{theorem}\label{extn}
The point processes $D_n,T_1,\ldots,T_n$ are independent
Poisson processes with respective intensities
$\lambda, \mu_1, \ldots, \mu_n$.
\end{theorem}
{\bf Proof.}
By Theorem~\ref{basic}, $D_1,T_1$ and $S_2$
are independent Poisson processes with respective intensities
$\lambda,\mu_1$ and $\mu_2$.  Applying Theorem~\ref{basic} again we see
that $D_2$ and $T_2$ are independent Poisson processes
with respective intensities $\lambda$ and $\mu_2$, and since $D_2$ and $T_2$
are determined by $D_1$ and $S_2$ they are independent of $T_1$.
Thus $D_2,T_1,T_2$ and $S_3$ are independent Poisson processes
with respective intensities $\lambda,\mu_1,\mu_2$ and $\mu_3$.  And so on.
The condition $\lambda < \min_{i\le n} \mu_i$ ensures that this procedure
is well-defined.
\hfill $\Box$

By repeated iteration of (\ref{pastk}) and (\ref{dmassk}), we obtain
(almost surely)
\begin{equation}\label{sum1}
Q_1(0)+\cdots +Q_n(0) = \sup_{s\ge 0} [\bar A(s)-(\bar S_1\i-c\cdots\i-c
\bar S_n)(s)].
\end{equation}

{\bf Remark.}  This formula, and variants of it, has been known for some
time now.
The first formula of this kind was observed by~\cite{muth} and later extended
in~\cite{sk,gop}.  Formulas of this kind can be used to compute first order
asymptotics for directed percolation (and polymer) variables, following a
program
introduced by Sepp\"al\"ainen~\cite{sepp} (see~\cite{noc} for a survey).
See also~\cite{hoc,oy1} for applications in a Brownian context.  We give a
brief
outline of these ideas in Section~\ref{perc}.

Iterating (\ref{dmassk}) we obtain, for each $k\le n$,
\begin{equation}\label{mass}
D_k(t)+Q_1(t)+\cdots +Q_k(t)=A(t)+Q_1(0)+\cdots +Q_k(0).
\end{equation}
We also have, by~(\ref{future}),
\begin{equation}
Q_k(t)= \sup_{u>t} [D_k(t,u)-T_k(t,u)] .
\end{equation}
Applying this repeatedly we obtain
\begin{equation}\label{sum}
Q_1(0)+\cdots +Q_n(0) = \sup_{t>0} [D_n(t)-(T_1\i-c\cdots\i-c T_n)(t)].
\end{equation}
Note that, on $\{ Q_1(0)+\cdots +Q_n(0)=0\}$,
\begin{equation}\label{star}
D_n(t)=(A\i-c S_1\i-c\cdots\i-c S_n)(t) ,
\end{equation}
and
\begin{equation}\label{star1}
T_k(t)=(S_k\s-c (A\i-c S_1\i-c\cdots\i-c S_{k-1}))(t),
\end{equation}
for $t\ge 0$, $k\le n$.

\section{Representations for non-colliding processes}\label{ncp}

We will start by showing how the two-dimensional analogue of
Theorem \ref{main} for Poisson processes is an immediate consequence
of Theorem~\ref{basic}, the extension of Burke's theorem given in
the previous section, and the symmetry formula~(\ref{future}).
As before, let $A$ and $S$ be Poisson processes with
respective intensities $\lambda$ and $\mu$.
\begin{theorem}\label{n=2}
The conditional law of $\{ (A(t),S(t)),\ t\ge 0\}$
given that $A(t)\le S(t)$ for all $t\ge 0$ is the same as the
unconditional law of $\{ \Gamma_2(A,S)(t),\ t\ge 0\}$.
\end{theorem}
{\bf Proof.}
By Theorem~\ref{basic} and the formula~(\ref{future}),
the conditional law of $\{ (A(t),S(t)),\ t\ge 0\}$
given that $A(t)\le S(t)$ for all $t\ge 0$ is the same as the conditional
law of $\{ (D(t),T(t)),\ t\ge 0\}$ given that $Q(0)=0$.
But when $Q(0)=0$, $(D(t),T(t))=\Gamma_2(A,S)(t)$ for $t\ge 0$.
Moreover, by~\eqref{past} and the independence of increments of
$A$ and $S$, $\{ \Gamma_2(A,S)(t),\ t\ge 0\}$ is independent of $Q(0)$.
Therefore, the conditional law of $\{ (A(t),S(t)),\ t\ge 0\}$
given that $A(t)\le S(t)$ for all $t\ge 0$ is the same as the
unconditional law of $\{ \Gamma_2(A,S)(t),\ t\ge 0\}$, as required.
\hfill $\Box$

Note that we can immediately deduce a discrete version of Pitman's
theorem, as follows.  If we set $X=S-A$, and $M(t)=\max_{0\le s\le t}X(s)$,
then $X$ is a simple random walk with positive drift, and the law of
$$2M-X = \Gamma_2(A,S)_2 - \Gamma_2(A,S)_1$$
is the same as that of $X$ conditioned to stay non-negative.

Now let $N^{(\mu)}=(N^{(\mu_1)}_1,\ldots,N^{(\mu_n)}_n)$ be the counting
functions of $n$
independent Poisson processes on $\R_+$ with respective intensities
$\mu_1<\mu_2<\cdots <\mu_n$.
That is, $N^{(\mu_k)}_k(t)$ is the measure induced by the $k^{th}$ Poisson
process on the interval $(0,t]$, with the convention that $N^{(\mu_k)}_k(0)=0$.
The following result was obtained in~\cite{oy2}.
\begin{theorem}\label{main1}
The conditional law of $N^{(\mu)}$, given that
$$N^{(\mu_1)}_1(t)\le\cdots\le N^{(\mu_n)}_n(t),\ \mbox{ for all }t\ge 0,$$
is the same as the unconditional law of $\Gamma_n(N^{(\mu)})$.
\end{theorem}
{\bf Proof.}  We prove this by induction on $n$;
the proof for $n=2$ is given above.
Assume that Theorem~\ref{basic} is true as stated
for a particular value of $n$, and moreover holds for any
choice of $\mu_1< \ldots < \mu_n$.  In the above setting we
have, by Theorem~\ref{extn}, that $D_n,T_1,\ldots,T_n$ are
independent Poisson processes with respective intensities
$\lambda, \mu_1, \ldots, \mu_n$.
Assume that $\lambda<\mu_1< \ldots < \mu_n$.
By the induction hypothesis, the conditional law of
\begin{equation} \{ (D_n(t),T_1(t),\ldots,T_n(t)),\ t\ge 0\} ,
\end{equation}
given that $T_1(t)\le\cdots\le T_n(t)$ for all $t\ge 0$,
is the same as the (unconditional) law of
\begin{equation} \{ (D_n(t),\Gamma_n(T_1,\ldots,T_n)(t)),\ t\ge 0 \} ;
\end{equation}
therefore, the conditional law of
\begin{equation} \{ (D_n(t),T_1(t),\ldots,T_n(t)),\ t\ge 0\} ,
\end{equation}
given that $D_n(t)\le T_1(t)\le\cdots\le T_n(t)$ for all $t\ge 0$,
is the same as the conditional law  of
\begin{equation} \{ (D_n(t),\Gamma_n(T_1,\ldots,T_n)(t)),\ t\ge 0 \} ,
\end{equation}
given that $D_n(t) \le (T_1\i-c\cdots\i-c T_n)(t)$ for all $t\ge 0$.
But, by~(\ref{sum}), this is precisely the condition that
$Q_1(0)+\cdots +Q_n(0)=0$ or, equivalently,
$Q_1(0)=\cdots =Q_n(0)=0$, and in this case we have, by~\eqref{star} and
\eqref{star1},
\begin{equation} (D_n(t),\Gamma_n(T_1,\ldots,T_n)(t))
=\Gamma_{n+1}(A,S_1,\ldots,S_n)(t)
\end{equation}
for $t\ge 0$; since this latter expression, by independence
of increments, is independent of $Q_1(0)+\cdots +Q_n(0)$, we are done.
\hfill $\Box$

Now let $N=(N_1,\ldots,N_n)$ be a collection of independent unit-rate
Poisson processes,
with $N(0)=(0,\ldots,0)$.  The function $h$ given by~(\ref{h}) is a
strictly positive harmonic
function for the restriction of the transition kernel of $N$ to the
discrete Weyl chamber
$E=W\cap\Z^n$.  This is well-known to harmonic analysists;
a proof is given in~\cite{kor}.  (In fact,
$h$ is harmonic for any random walk, or any
L\'evy process, with exchangeable increments, provided the required moments
exist.)
Let $x^*+\hat N$ be a realisation of the corresponding Doob $h$-transform,
started at
$x^*=(0,1,\ldots,n-1)\in E$ (so that $\hat N(0)=(0,\ldots,0)$).
Apart from providing a convenient framework in which to apply Donsker's theorem
and deduce Theorem~\ref{main} from Theorem~\ref{main1}, the process $\hat
N$ is interesting
in its own right.
In~\cite{kor} it is shown that the random vector $\hat N(1)$ is distributed
according to the
{\em Charlier} ensemble, a discrete orthogonal polynomial ensemble.
That is,
\begin{equation}\label{ce}
P(\hat N(t)=y)=C_t h(x^*+y)^2 P(N(t)=y),
\end{equation}
where $C_t$ is the same (!) normalisation constant as in the Brownian case,
given by~(\ref{ct}).
Thus, the next result, which follows from Theorem~\ref{main1}, yields
a representation for the Charlier ensemble.  For more on discrete
orthogonal polynomial ensembles, see~\cite{j2}.

\begin{theorem}\label{main2}
The processes $\hat N$ and $\Gamma_n(N)$ have the same law.
\end{theorem}

Some elementary potential theory is needed to deduce this from
Theorem~\ref{main1},
since the function $h$ is not the only positive harmonic function for
the restriction of the transition kernel of $N$ to the discrete Weyl chamber
$E=W\cap\Z^n$.  We refer the reader to~\cite{oy2} for details of the proof,
and to~\cite{kor} for the necessary asymptotic analysis of the Green's
function.
We remark that, for any $y\in\R_+^n$, the function defined by
$g_y(x)=\mbox{Schur}_x(y)$
is a positive harmonic function for $N$ on $E$, and note that $g_y$ is a
multiple
of $h$ if all the components of $y$ are equal.  This follows from the analysis
of the Green's function presented in~\cite{kor} (alternatively it can be
deduced from
more general representation theoretic arguments, as in~\cite{biane2}).

Theorem~\ref{main} (the Brownian version)
follows from Theorem~\ref{main2} by careful application
of Donsker's theorem (see~\cite{oy2} for details).

The analogues of Theorems~\ref{main1} and \ref{main2} for discrete-time
walks are presented
in~\cite{kor}; in this case, replace `Poisson' by `binomial', `Charlier' by
`Krawtchouk',
and modify the definition of $\Gamma_n$.  The proofs are similar, but there are
additional complications due to the fact that the walkers can jump
simultaneously;
the asymptotic analysis of the Green's function is also more difficult in
this case.

\section{Output theorems generally}\label{output}

Theorem~\ref{basic} seems to be a special case of a more general
phenomenon.  In this section we will present several examples
of `output theorems'; in all cases, reversibility and symmetry
play a key role.

The analogue of Theorem~\ref{basic} holds for Brownian motions with
drift.  Suppose now that $A$ and $S$ are independent Brownian motions,
indexed by $\R$, with respective drifts $\lambda <\mu$. That is, for each $s$,
$\{A(s+t)-A(s),\ t\ge 0\}$ is a Brownian motion with drift $\lambda$.
Set $A(s,t)=A(t)-A(s)$, $S(s,t)=S(t)-S(s)$, $$Q(t) = \sup_{s\le t}
[A(s,t)-S(s,t)],$$
$$D(t) = A(t)+Q(0)-Q(t),$$ $$T(t) = S(t)-Q(0)+Q(t),$$ for $t\in\R$.  Then
the pair
$(D,T)$ has the same law as $(A,S)$.  Several proofs of this fact are given
in~\cite{oy1};
an equivalent result was earlier presented in~\cite{hw} (see below).

Note that Theorem~\ref{basic} can also be stated as follows, by considering
the process $X=S-A$ and its transform $\tilde X=T-D$.
\begin{theorem}\label{basic-srw}
Let $X$ be a simple random walk (in continuous time) with positive drift,
indexed by $\R$, and set $$Q(t) = \sup_{-\infty <s< t} X(s)-X(t),$$
$$\tilde X(t) =  X(t) + 2[Q(t)-Q(0)].$$  Then $\tilde X$ has the same law
as $X$.
\end{theorem}
Note that, in this setting, the `symmetry formula' \eqref{future} becomes
\begin{equation}\label{symm}
Q(t)=\sup_{u>t} \tilde X(t)-\tilde X(u).
\end{equation}
Again, the Brownian version of this result holds: this was
the Brownian version of Burke's theorem presented in~\cite{hw}.

The analogue of Theorem~\ref{extn} can also be shown to hold for
Brownian motions with drifts, by exactly the same argument given in the
proof of Theorem~\ref{extn}.

We will now describe some output theorems, which are completely
analogous to those described above, but in the usual algebra
(as opposed to the max-plus algebra).  These were given in~\cite{oy1},
and motivated by the recent extensions of Pitman's $2M-X$ theorem
obtained by Matsumoto and Yor~\cite{my} in the context of exponential
functionals of Brownian motion.  In words, the analogues of Theorems
\ref{basic} and \ref{basic-srw} hold for Brownian motions if we replace
`$\sup$' by `$\log\int\exp$' in the definition of $Q$.  More precisely,
suppose that $A$ and $S$ are Brownian motions with respective drifts
$\lambda<\mu$ and we set, for $t\in\R$,
$$Q(t) = \log\int_{-\infty}^t \exp [A(s,t)-S(s,t)] ds$$
$$D(t) = A(t)+Q(0)-Q(t),$$ $$T(t) = S(t)-Q(0)+Q(t).$$ Then the pair
$(D,T)$ has the same law as $(A,S)$.
As discussed in~\cite{oy1}, symmetry and reversibility
play a big role here.  In particular, the process $Q$ is reversible, and
the analogue of the symmetry formula \eqref{symm} holds:
\begin{equation}\label{symm-le}
Q(t)=\log\int_t^\infty \exp [D(t,u)-T(t,u)] du.
\end{equation}
Let us return to the max-plus version for a moment:
intuitively, the statement of Theorem~\ref{basic} can almost
be seen as a direct consequence of the reversibility of the process $Q$
combined with the symmetry formula \eqref{future};  the only thing
which makes this difficult is the non-invertibility of the transformation.
This problem does not exist in the usual algebra, and in fact the statement
that $(D,T)$ has the same law as $(A,S)$ (when everything is defined using
$\log\int\exp$) is actually equivalent to the statement that $Q$ is reversible.
For a more precise statement of this equivalence, which holds quite generally,
see~\cite{oy1}.  The version of Pitman's $2M-X$ theorem obtained in this
context
by Matsumoto and Yor~\cite{my} states that, if $X$ is a Brownian motion
with drift and
$$2M(t) = \log\int_0^t e^{2 X(s)} ds,$$
then $2M-X$ is a certain diffusion process.
Note that by following the arguments given at the start of the last section
to deduce Pitman's theorem for random walks from the output
theorem~\ref{basic},
we can (formally) recover the Matsumoto-Yor representation and, moreover, see
that, in some suitable sense, the law of the diffusion $2M-X$ can be
interpreted as
the conditional law of $X$, given that $$\int_0^\infty e^{2 X(s)} ds = 1.$$

Similarly, the analogue of Theorem~\ref{extn} can be shown to hold in the
$\log\int\exp$-world for Brownian motions with drifts, by exactly the same
extension argument given in the previous section.  Note also that the arguments
given in the proof of Theorem~\ref{main1} can also be (formally) carried over
to obtain a $\log\int\exp$-analogue of Theorem~\ref{main}.  This has
the following form:  define a mapping $\Pi_n$ analogously with $\Gamma_n$
but replacing $\sup$ by $\log\int\exp$ (and $\inf$ by $-\log\int\exp(-)$),
in the definitions of $\i-c$ and $\s-c$.  Then, by these formal arguments,
the law of $\Pi_n(B)$ can be interpreted as the conditional law of $B$, given
that
$$\int_0^\infty \exp[B_i(s)-B_{i+1}(s)] ds = 1,$$
for each $i=1,\ldots,n-1$.

In~\cite{hmo}, the following {\em non-Markovian} version of
Pitman's theorem is obtained using similar ideas:  reversibility arguments do
not require the Markov property.

\begin{theorem}
Let $(\xi_k,\ k\ge 0)$ be a Markov chain on $\{-1,+1\}$ with $\xi_0=1$ and
transition probabilities
$P(\xi_{k+1}=1|\ \xi_k=1)=a$ and \hbox{$P(\xi_{k+1}=-1|\ \xi_k=-1)=b<a$.}
Set $X_0=0$, $X_n=\xi_1+\cdots +\xi_n$ and $M_n=\max_{0\le k\le n}X_k$.
The process $2M-X$ has the same law as that of $X$ conditioned to stay
non-negative.
\end{theorem}
The continuous analogue of the above theorem, with $\xi$ replaced by a zig-zag
process, is also given in~\cite{hmo}.

We end this section with one final example, quite different from the above,
to give an indication of how generally these output theorems, and the
corresponding representation theorems, hold.
It is related to the first-order autoregressive process.
Let $(b_n,\ n\in\Z )$ be a sequence of independent standard normal
random variables, and fix $0<a<1$.
Consider the stationary process $(X_n,\ n\in\Z )$ defined by
$X_n=\sum_{j=0}^\infty a^{-j} b_{n-j}$.  Then $X_{n+1}=aX_n+b_n$ for all $n$.
Now observe that, if we set
$$\hat b_n = X_n - aX_{n+1} ,$$
then $X_n=aX_{n+1}+\hat b_n$ and $X_n=\sum_{j=0}^\infty a^{-j} \hat b_{n+j}$.
Now, by the symmetry of this construction, and the fact that the process
$X$ is reversible, we see that:
\begin{theorem} The sequences $b$ and $\hat b$ have the same law; that is,
$(\hat b_n,\ n\in\Z )$ is a sequence of independent standard normal random
variables.
\end{theorem}
This is the output theorem.  For the analogue of
Pitman's theorem in this case, set $Y^{(x)}_0=x$ and, for $n\ge 0$,
$Y^{(x)}_{n+1}=aY^{(x)}_n+b_n$,
$$d^{(x)}_n = Y^{(x)}_n - aY^{(x)}_{n+1} .$$  Consider the conditional law
of $(b_n,\ n\ge 0)$ given that $\sum_{n=0}^\infty a^n b_n =x$.
(It is easy to see that this is well-defined.)  From the output theorem,
and the symmetry formula $X_0=\sum_{n=0}^\infty a^{-n} \hat b_n$,
this is the same as the law of $(\hat b_n,\ n\ge 0)$ given that $X_0=x$.
But, when $X_0=x$, $\hat b_n = d^{(x)}_n$ and, moreover, $(d^{(x)}_n,\ n
\ge 0)$
is independent of $X_0$; it follows that
\begin{theorem}
The law of the sequence $(d^{(x)}_n,\ n \ge 0)$ is the same as
the conditional law of
$(b_n,\ {n \ge 0})$ given that $\sum_{n=0}^\infty a^n b_n =x$.
\end{theorem}

\section{First-order asymptotics for directed percolation and directed
polymers}
\label{perc}

The formula
\begin{equation}\label{sum2}
Q_1(0)+\cdots +Q_n(0) = \sup_{s\ge 0} [\bar A(s)-(\bar S_1\i-c\cdots\i-c
\bar S_n)(s)],
\end{equation}
introduced in Section~\ref{burke}, is very useful for computing first-order
asymptotic results for directed percolation variables.

In Section~\ref{burke}, $A,S_1,\ldots, S_n$ are independent
Poisson processes with respective rates $\lambda,\mu_1,\ldots,\mu_n$, with
$\lambda < \min_i\mu_i$.
In this case, the random variables $Q_i(0)$ are independent and
geometrically distributed with
respective parameters $\lambda/\mu_i$.  Assume $\lambda<1$.
If we let $\mu_i=1$, for each $i$, and divide (\ref{sum2}) by $n$, we can
(in principle)
let $n\to\infty$ to obtain:
\begin{equation}\label{vf}
\frac{\lambda}{1-\lambda} = \sup_{x>0} [\lambda x-\gamma(x)],
\end{equation}
where, almost surely,
\begin{equation}
\gamma(x)=\lim_{n\to\infty} \frac{1}{n}(S_1\i-c\cdots\i-c S_n)(xn).
\end{equation}
(Here we have implicitly used the reversibility of $A$, and of $S_i$, for
each $i$.)
The existence of $\gamma$ follows from Kingman's subadditive ergodic theorem.
But (\ref{vf}) is essentially a Legendre transform.  In particular, it can
be inverted to obtain $$\gamma(x)=(\sqrt{x}-1)^2 1_{x>1}.$$

This ingenuous idea is due to Sepp\"al\"ainen~\cite{sepp}, who used it to
compute asymptotics for a certain last-passage Bernoulli directed
percolation problem.
See~\cite{noc} for a survey on the application of this technique to a variety
of directed percolation problems.  In each context, there is an associated
`growth model', and these limit theorems can also be interpreted as determining
the limiting shape of the associated growth model.  See also~\cite{j1} for
related work in this area.

For the Brownian model, we have, almost surely,
\begin{equation}\label{zoot}
\lim_{n\to\infty} \frac{1}{n}(B_1\i-c\cdots\i-c B_n)(xn)=2\sqrt{x},
\end{equation}
where $B_1,B_2,\ldots$ are independent standard Brownian motions.
This is described in~\cite{oy1} and proved with sharp uniform concentration
estimates in~\cite{hoc}.  The corresponding corner growth model is also
discussed in that paper.  In this case, the formula (\ref{sum2}) holds
with $A,S_1,\ldots, S_n$ replaced by independent standard Brownian motions
with respective drifts $\lambda,1,\ldots,1$ and the $Q_i(0)$ are i.i.d.
exponentially
distributed with parameter $\lambda$.
Note that the limiting result (\ref{zoot}) can also be deduced from
Theorem~\ref{bgtw}
by applying standard results from random matrix theory.

The formula (\ref{sum2}) is also valid when `$\sup$' is replaced by
`$\log\int\exp$'.  In this context it is used in~\cite{oy1} to compute
(formally) the free energy density associated with a certain directed polymer
in a random environment.  More precisely, for $\beta>0$, if
$$Z_n(\beta)=\int_{-n<s_1<\cdots <s_{n-1}<0} ds_1\ldots ds_{n-1}
\exp\left\{ \beta (B_1(-n,s_1)+\cdots +B_n(s_{n-1},0))\right\},$$
where $B_i(s,t)=B_i(t)-B_i(s)$, then, almost surely,
$$f(\beta) =\lim_{n\to\infty}\frac{1}{n}\log Z_n(\beta)
= -g(-\beta^2)-2\log\beta ,$$
where
$$g(x)=\sup_{y>0} [xy+\Psi(y)],$$
for $x<0$, and $\Psi=\Gamma'/\Gamma$ is the digamma function.

{\em Acknowledgements.}  I wish to thank everyone who has participated
in this research and made it so enjoyable.
This article was completed during my stay at the Laboratoire de Probabilit\'es,
Paris 6.  I would like to thank everyone at the Laboratoire, and Marc Yor in
particular, for their kind hospitality.  Thanks also to the CNRS for financial
support.

\end{document}